\documentclass[12pt,reqno,sumlimits]{amsart}   
\usepackage{amssymb,amscd,amsmath,epsfig}
\usepackage{%
  array,
  booktabs,
  dcolumn,
  rotating,
  shortvrb,
  tabularx,
  units,
  url,
}
\newtheorem{thm}{Theorem}[section]
\newtheorem{cor}[thm]{Corollary}
\newtheorem{lem}[thm]{Lemma}
\newtheorem{prop}[thm]{Proposition}
\theoremstyle{definition}
\newtheorem{defn}{Definition}[section]
\newtheorem{rem}{Remark}[section]

\newcommand{\Pdo}{{\Psi {\rm DO}}}

\newcommand{\R}{{\mathbb R}}
\newcommand{\G}{{\rm Ell}_0^*}

\newcommand{\C}{{\mathbb C}}

\newcommand{\Z}{{\mathbb Z}}

\newcommand{\g}{{\mathfrak  g}}

\newcommand{\calB}{{\mathcal B}}

\newcommand{\calE}{{\mathcal E}}

\newcommand{\calG}{{\mathcal G}}

\newcommand{\calR}{{\mathcal R}}

\renewcommand{\to}{\longrightarrow}

\newcommand{\ev}{\operatorname{ev}}

\newcommand{\Tr}{\operatorname{Tr}}
\newcommand{\tr}{\operatorname{Tr}}

\newcommand{\con}[2]{\nabla^{#1}_{#2}}

\newsavebox{\savepar}

\newcommand{\ip}[1]{\langle #1 \rangle}

\numberwithin{equation}{section}

\newcounter{labelflag} 
\newcommand{\labelon}{\setcounter{labelflag}{1}}
\newcommand{\Label}[1]{
                       \ifnum\thelabelflag=1
                          \ifmmode
                             \makebox[0in][l]{\qquad\fbox{\rm#1}}
                          \else
                             \marginpar{\vspace{0.7\baselineskip}
                                        \hspace{-1.1\textwidth}
                                        \fbox{\rm#1}}
                          \fi
                       \fi
                       \label{#1}
                      }
\labelon

\begin{document}
\title{Infinite Dimensional Chern-Simons Theory}
\author[S. Rosenberg]{Steven Rosenberg}
\address{Department of Mathematics and Statistics\\
  Boston University}
\email{sr@math.bu.edu}
\author[F. Torres-Ardila]{Fabi\'an Torres-Ardila}
\address{Department of Mathematics and Statistics\\
  Boston University}
\email{fatorres@math.bu.edu}

\begin{abstract}  We extend  finite
dimensional Chern-Simons theory to certain infinite dimensional principal
bundles with connections, in particular to 
the frame bundle $FLM\to LM$ over the
loop space of a Riemannian manifold $M$. Chern-Simons forms are defined
roughly as
in finite dimensions with the invariant polynomials replaced by
appropriate Wodzicki residues. This produces odd dimensional $\R/\Z$-valued
cohomology classes on $LM$ if $M$ is parallelizable.  
We compute an example of a metric on the loop
space of $S^3\times S^1$ for which the three dimensional Chern-Simons class is
nontrivial.
\end{abstract}
\maketitle


     \section{{\bf Introduction}}

A principal $G$-bundle $E\to B$
with connection has characteristic classes determined
by applying invariant polynomials $P$ 
on the Lie algebra of $G$ to the curvature $\Omega$
of the connection.  The apparently worst possible case, when
the form $P(\Omega)$ vanishes pointwise, in fact leads to nontrivial
secondary classes, the Chern-Simons classes, in $H^{\rm odd}(B,\R/\Z)$
\cite{C-S}.  For the frame bundle of a Riemannian manifold,
these classes contain geometric information.

This theory assumes that $G$ is a finite dimensional Lie group.  For
interesting infinite dimensional Riemannian manifolds such as spaces of
maps $C^\infty(N,M)$ between Riemannian manifolds, 
and in particular for loop spaces $LM$,
the frame bundle has as structure group the gauge group $\calG$
of an auxiliary finite rank bundle.  However, the Levi-Civita connection on
mapping spaces takes values in pseudodifferential operators ($\Pdo$s) on the
auxiliary bundle, so it is necessary to expand the structure group to $\G$, 
the group of zeroth order invertible $\Pdo$s.

In \cite{P-R2}, a theory of characteristic classes was developed for $\calG$-
and $\G$-bundles.
The invariant
polynomials for classical groups are generated by $A\mapsto \tr(A^l)$, but
the 
$\Pdo$s for mapping spaces are often not trace class operators.  Thus the
invariant polynomials in \cite{P-R2} use alternate traces, namely
the leading order symbol or the Wodzicki residue, to construct characteristic
forms and classes.  While there are examples of nontrivial
leading order Chern classes, no nontrivial examples of the more natural
Wodzicki-Chern classes are known.

Both the leading order Chern forms (more precisely, Pontrjagin forms)
and the Wodzicki-Chern forms 
vanish on the frame bundle $FLM\to LM$ for loop
spaces (except possibly the first Wodzicki-Chern form -- see \S4),
so it is natural to
look for Chern-Simons classes.  This paper develops a theory of Chern-Simons
classes for loop spaces for the two types of characteristic classes.  For the
frame bundle, we are unable to produce nonzero 
leading order Chern-Simons classes.
However, the first Wodzicki-Chern-Simons class is nontrivial in 
$H^3(L(S^3\times S^1), \R/\Z)$ for a wide class of metrics, 
which shows that this theory is nonvacuous.  

The paper is organized as follows.  In \S2, relevant results on the Riemannian
geometry of $FLM$ are collected.  
In \S3, finite dimensional Chern-Simons theory is reviewed, and the theory is
extended to infinite dimensional bundles with structure group either 
$\calG$ or $\G$.
\S4 gives an integrality result for leading order Chern classes (this is joint
work with S. Paycha) and new vanishing results for these classes.  The
integrality result and the existence of a classifying space with universal
connection for the gauge group are crucial to define Chern-Simons classes for 
$\calG$-bundles.  Unfortunately, we do not know if these
results hold
for $\G$-bundles.  As a result, we can only define 
Chern-Simons classes for trivial $\G$-bundles,
such as loop spaces of parallelizable manifolds, and the
Chern-Simons class apparently depends on the choice of a global frame.  
This restriction to parallelizable manifolds often occurs even
in finite
dimensions, as Chern-Simons classes are notoriously difficult to compute.
In any
case, in \S5 a strategy for proving the nontriviality of a Chern-Simons class
is given, and in \S6 this is applied to produce a nontrivial 
Wodzicki-Chern-Simons class on $L(S^3\times S^1).$

Helpful conversations with David Fried and Sylvie Paycha are gratefully
acknowledged, particularly in \S4.1.

\section{{\bf Preliminaries on the Geometry of $LM$}}


Let $(M, \langle\ ,\ \rangle)$ 
be a compact oriented Riemannian $n$-manifold with  loop space $LM
= C^\infty(S^1,M)$ of smooth loops. 
$LM$ is a smooth infinite dimensional manifold, but it is technically simpler 
to work with the smooth manifold of loops in some Sobolev class $s \gg 0,$
as we now recall. For $\gamma\in LM$, the formal
tangent space $T_\gamma LM$ is 
$\Gamma(\gamma^*TM)$, the space
 of smooth sections of the pullback bundle $\gamma^*TM\to
S^1$. For $s>1/2$, we complete $\Gamma(\gamma^*TM)$
 with respect to the Sobolev inner product
 \begin{equation*}\label{eq:Sob1}
\langle X,Y\rangle_{s}=\frac{1}{2\pi}\int \langle(1+\Delta)^s
X(\alpha),Y(\alpha)
\rangle_{\gamma (\alpha)}d\alpha,\  X,Y\in \Gamma(\gamma^*TM).
\end{equation*}
Here $\Delta=D^*D$, with $D=D/d\gamma$ the covariant derivative along
$\gamma$. We denote this completion by  $H^s(\gamma^*TM)$. 

A small neighborhood $U_\gamma$
of the zero section in $H^s(\gamma^*TM)$ is a
coordinate chart near $\gamma$ in the space of $H^s$ loops
via the pointwise exponential map
\begin{equation*}
\begin{split}
\exp_\gamma&:U_\gamma
\to L M, \ X \mapsto 
\left(\alpha\mapsto \exp_{\gamma(\alpha)} X(\alpha)\right).  
\end{split}
\end{equation*}
It is easy to check that $\exp_\gamma$ is a bijection onto its image.
The differentiability of the transition functions $\exp_{\gamma_1}^{-1}\cdot
\exp_{\gamma_2}$ is proved in
\cite{E} and \cite[Appendix A]{Freed1}.
Since 
$\gamma^*TM$
is (non-canonically) isomorphic to the trivial bundle ${\mathcal R}^n =
S^1\times \R^n\to S^1$, 
the model space for $LM$ is the set of 
$H^s$ sections of this trivial bundle.  

\subsection{The Tangent Bundle TLM}

The tangent bundle
$TLM$ has
transition functions $d(\exp_{\gamma_1}^{-1} 
\circ \exp_{\gamma_2})$.  Under
the isomorphisms $T_{\gamma_1}LM \simeq {\mathcal R}^n \simeq
T_{\gamma_2}LM$, the transition functions are gauge transformations of
${\mathcal R}^n.$

The $H^s$
metric makes $LM$ a Riemannian manifold. 
The $H^s$
Levi-Civita connection on $LM$ is determined by the six term formula
\begin{eqnarray*}\label{eq:SixTerm}
\ip{\con{s}{Y}X,Z}_s &=& X\ip{Y,Z}_s+Y\ip{X,Z}_s-Z\ip{X,Y}_s\\
&&\qquad +\ip{[X,Y],Z}_s+\ip{[Z,X],Y}_s-\ip{[Y,Z],X}_s.\nonumber
\end{eqnarray*}
$\con{s}{}$ is explicitly computed in \cite[Theorem 2.2]{M-R}:
\begin{eqnarray}\label{eq:sconn}
2(\nabla^s_XY)^a &=&
 \ 2(\delta_X Y)^a+(1+\Delta)^{-s}\left[ g^{af}\delta_Y g_{ef}
  ((1+\Delta )^s Y)^e +(\delta_X (1+\Delta)^s Y)^a\right]\nonumber\\
&&\qquad +(1+\Delta)^{-s}\left[ g^{af}\delta_Y g_{bf}
  ((1+\Delta )^s X)^b +(\delta_Y(1+\Delta)^s X)^a\right]\\
&&\qquad -(1+\Delta)^{-s}\left[ g^{at}\delta_t g_{bf}
  ((1+\Delta )^s X)^bY^f +g^{at} g_{bf}((1+\Delta)^s X)^bY^f\right],\nonumber
\end{eqnarray} 
where $\delta_X$ is the variation in the $X$ direction, and the
$g$'s and $\Gamma$'s are the metric tensor and Christoffel symbols of 
 $M$ in a local frame
$\{e_a(\gamma_1,\alpha)\in T_{\gamma_1(\alpha)}M\}$ for $\gamma_1$ in a
  neighborhood of $\gamma$. 
By
(\ref{eq:sconn}), 
the connection one-form of the
operator $Y\mapsto \nabla_XY$ is a pseudodifferential operator 
($\Pdo$) of order $0$ acting on sections of $\gamma^*TM$. The curvature
$\Omega^{s}(X,Y)=\con{s}{X}\con{s}{Y}-\con{s}{Y}\con{s}{X}-\con{s}{[X,Y]}$,
as an operator $Z\mapsto \Omega^{s}(X,Y)Z$, is a  $\Pdo$ of order -1
 \cite{M-R}.

\subsection{The Frame Bundle $FLM$}

The frame bundle $FLM\to LM$ is constructed
as in the finite dimensional case. The
fiber over $\gamma$ is isomorphic to the gauge group $\calG$ of $\calR^n$
and fibers are glued by the transition functions for
$TLM$. Thus the frame bundle is
topologically a
$\calG$-bundle.

However, the connection and curvature forms take values in $\Pdo_{\le 0}$, 
the set of $\Pdo$s of order at most zero. These
forms should take values in the Lie algebra of the structure
group, so we consider the extended structure group $\G$, the
group of classical invertible (and therefore elliptic)
$\Pdo$s of order 0 acting on $\calR^n$,
as $\Pdo_{\le 0}=\operatorname{Lie }(\G) $. Note that
$\calG$ embeds in $\G$ as multiplication
operators. Thus $(FLM,\theta^s)$ as a geometric
bundle (i.e.~as a  bundle with connection $\theta^s$ associated to
$\nabla^s$, as explained below) is a $\G$-bundle.

In summary, we have
$$ \begin{array}{ccc}
\calG&\longrightarrow &FLM\\
& & \downarrow\\
& & LM
\end{array}
\ \ \ \ \ \ \ \ \ 
\begin{array}{ccc}
\G&\longrightarrow &(FLM,\theta^s)\\
& & \downarrow\\
& & LM.
\end{array}
$$

\subsection{Connections on the Frame Bundle}

We summarize the relationship between 
 the Levi-Civita connection $\theta^s$ on the frame bundle and
local expressions for the Levi-Civita connection on the tangent
bundle. Let $\chi:N\to FN$ be a local section of the frame bundle of an
$n$-manifold $N$. A metric  connection $\nabla$ on $TN$ with local
connection one-form $\omega$ determines a connection $\theta_{FN}\in
 \Lambda^1(FN, {\mathfrak so}(n))$ on $FN$ 
by {\it (i)} 
the standard property for $\theta_{FN}$ on vertical vectors, and {\it (ii) }
$\theta_{FN}(Y_u)=\omega (X_p),$ for $ Y_u=\chi_*X_p$
\cite{Spi}, or equivalently
\begin{equation}\label{localexp}
\chi^*\theta_{FN} = \omega.
\end{equation}
For $N=LM$, $\nabla^s$ determines a connection $\theta^s\in
\Lambda^1(LM, \Pdo_{\leq 0})$ on $FLM$, and the computations of the symbols
$\sigma_0(\theta^s)$, $\sigma_{-1}(\theta^s)$ reduce to the computations 
of 
$\sigma_0(\omega^s)$ and $\sigma_{-1}(\omega^s)$ of the
Levi-Civita connection one-form. 
By \cite{M-R}, for $X\in \Gamma(
\gamma^*TM)$, we have
\begin{equation}\label{eq:zero}
\sigma_0(\omega^s(X))^k_l
= \frac{1}{2}\left(\Gamma^{k}_{lp}
+g^{kb}g_{lm}\Gamma^{m}_{bp}\right)X^p,
\end{equation}
\begin{eqnarray}\label{eq:Minus}
\sigma_{-1}(\omega^s(X))^k_l &=& 2is \xi^{-1}\left[ \partial_m
  \Gamma^k_{ln}\dot\gamma^n X^m
+\Gamma^k_{lm}\dot X^m+\Gamma^k_{ma}\Gamma^a_{ln}
  X^m\dot\gamma^n-\Gamma^k_{pn}\Gamma^p_{ml}X^m\dot\gamma^n \right.\nonumber\\
  &&\ \  \left.
 +\dot\gamma^n\partial_n\Gamma^k_{ml}X^m 
  + \Gamma^k_{ml}\dot
  X^m-\Gamma^k_{pn}g^{pb}\Gamma^r_{bt}
X^tg_{rl}\dot\gamma^n-
\Gamma^p_{mn}g_{lp}
\Gamma^m_{br}X^r\dot\gamma^n \right.\nonumber\\
&&\ \  \left.
+ \left(
  g^{kb}\Gamma^p_{bn}X^n\right)\dot{} g_{p
  l}\right],
\end{eqnarray}
where $\dot{} = d/d\alpha$ along $\gamma.$
Using a frame $\{\psi_a\}$ dual to  $\{e_a\}$, we may express 
e.g.~(\ref{eq:zero}) as a form:
\begin{equation}\label{eq:zero1}
\begin{split}
 \sigma_0(\omega^s)^k_l &=
\frac{1}{2}\left[\Gamma^k_{lp}+g^{kb}g_{ls}\Gamma^s_{bp} \right]\psi^p
\end{split}
\end{equation}

\section{{\bf Chern-Simons Forms}}

In this section, the essentials of Chern-Simons theory
are reviewed and extended to infinite dimensions. 
First, we recall the theory of characteristic classes in finite
dimensions and the main results of the seminal paper \cite{C-S}. 

\subsection{Characteristic Classes}

Consider a finite dimensional Lie group $G$, a manifold $M$, and a $G$-bundle
 $E\to M$.
For $\g={\rm Lie }(G)$, set $
 \g^l=\g^{\otimes l}$ and
\begin{equation*}I^l(G)
= \{P:\g^l\to \R\ | P\ \text{symmetric,
  multilinear, ad-invariant}\}.
\end{equation*}
For $\phi\in\Lambda^k(E,\g^l)$, $P\in I^l(G)$, set
$P(\phi)=P\circ \phi\in\Lambda^k(E)$. Two key 
properties are:
\begin{itemize}
\item (The {\it commutativity property}) 
For $\phi\in\Lambda^{k}(E,\g^l)$, 
\begin{equation}\label{eq:deri}
d(P(\phi))=P(d\phi).
\end{equation}
\item (The {\it infinitesimal invariance property})
For $\psi_i\in\Lambda^{k_i}(E,\g)$, $\phi\in\Lambda^{1}(E,\g)$ and $P\in
  I^l(G)$, 
\begin{equation}\label{eq:inva}
\sum^l_{i=1}
(-1)^{k_1+\dots+k_i}P(\psi_1\wedge\dots\wedge[\psi_i,\phi]\wedge\dots
\psi_l)=0. 
\end{equation}
\end{itemize}

\begin{rem}
For classical Lie groups $G$, (\ref{eq:deri}) and (\ref{eq:inva}) follow
trivially from
the fact that the polynomials $I^l(G)$ is generated by
the Newton polynomials $\Tr(A^l)$ and properties of the ordinary finite
dimensional trace.
\end{rem}

\begin{thm}[The Chern-Weyl Homomorphism \cite{K-N}] 
Let $E\to M$ have a connection $\theta$ with curvature $\Omega_E\in
\Lambda^2(E,\g)$. For $P\in I^l(G)$, $P(\Omega^l_E)$ is a closed
 invariant real form on $E$, and so
determines a closed form
$P(\Omega_M)\in \Lambda^{2l}(M,\R)$, a form with lift
$P(\Omega^l_E)$. 
The Chern-Weil map
\begin{equation*}
\oplus_{l=1}I^l(G)\to H^{*}(M,\R), \ P\mapsto [P(\Omega_M)]
\end{equation*}
is a well-defined algebra homomorphism.
\end{thm} 

$[P(\Omega_M)]$ is
called the {\it characteristic class} of $P$. 

\subsection{Chern-Simons Theory for Finite Dimensional Bundles}

A crucial observation of Chern-Simons is that $P(\Omega^l_E)$ is exact,
although in general  $P(\Omega_M)$ is not.

\begin{prop}\label{eq:CS}\cite[Proposition 3.2]{C-S} Let $G$ be a finite
  dimensional Lie group.
For a $G$-bundle $E\overset{\pi}{\to} M$ with connection $\theta$
and curvature $\Omega = \Omega_E$, and for $P\in I^l(G)$, set
\begin{equation*}\label{eq:ChernSimons}
\begin{split}
\phi_t&=t\Omega+\frac{1}{2}(t^2-t)[\theta,\theta]\\
TP(\theta)&=l\int_0^1 P(\theta\wedge \phi^{l-1}_t)dt.
\end{split}
\end{equation*}
Then $dTP(\theta)=P(\Omega)\in \Lambda^{2l}(E)$.
\end{prop}

\begin{proof}
We recall the proof for later purposes. 
Set $f(t)=P(\phi^l_t)$, so
$P(\Omega^l)=\int_0^1f^\prime(t)dt$. We show $f^\prime(t)=l\cdot
dP(\theta\wedge \phi^{l-1}_t)$ by computing each side. First, we have
\begin{equation}\label{eq:lsid}
\begin{split}
f^\prime(t)&=\frac{d}{dt}\left(P(\phi^l_t)\right)
=P\left(\frac{d}{dt}\phi^l_t\right)=lP\left(\frac{d}{dt}
\phi_t\wedge\phi^{l-1}\right)\\
&=lP(\Omega\wedge\phi^{l-1})+l\left(t-\frac{1}{2}\right)P\left([\theta,\theta]\wedge\phi^{l-1}_t\right),
\end{split}
\end{equation}
where we have used the commutativity property (\ref{eq:deri}). On the other
hand, we have
\begin{eqnarray}\label{eq:rsid}
l\cdot
dP(\theta\wedge \phi^{l-1}_t)&=& lP(d\theta\wedge 
\phi^{l-1}_t)-l(l-1)P(\theta\wedge d\phi_t\wedge \phi^{l-2}_t)\\
&=&  lP( \Omega\wedge \phi^{l-1}_t)-\frac{1}{2}l
P([\theta,\theta]
\wedge \phi^{l-1}_t)-l(l-1)P(\theta\wedge d\phi_t\wedge
\phi^{l-2}_t),\nonumber 
\end{eqnarray}
by (\ref{eq:deri}) and the structural equation
$\Omega = d\theta + \frac{1}{2}[\theta,\theta].$
 Since $d\phi_t=t[\phi_t,\theta]$, the last term
in (\ref{eq:rsid}) equals
\begin{equation*}
l(l-1)P(\theta\wedge d\phi_t\wedge \phi^{l-2}_t)
=l(l-1)P(\theta\wedge t[\phi_t,\theta]\wedge \phi^{l-2}_t).
\end{equation*}
Using the invariance property (\ref{eq:inva}) with $\phi=\theta$,
$\psi_1=\theta$ and $\psi_k=\psi_t, k=2,\dots, l-1$, we obtain
\begin{equation*}
l(l-1)P(\theta\wedge t[\phi_t,\theta]\wedge \phi^{l-2}_t)=-lt
P([\theta,\theta]\wedge\phi^{l-1}_t).
\end{equation*}
This  implies (\ref{eq:rsid}) equals (\ref{eq:lsid}).
\end{proof}

Setting $M=BG$ in the theorem gives the 
{\it universal Chern-Weil homomorphism}
\begin{equation*}
  W:I^l(G)\to H^{2l}(BG,\R).
  \end{equation*}
We write $P\in I^l_0(G)$ if $W(P)\in H^{2l}(BG, \Z)$. For this subalgebra of
polynomials, we obtain more information on $TP(\theta)$.

\begin{thm}\label{thm:Coh}\cite[Proposition 3.15]{C-S}. Let $E\to B$ be a 
$G$-bundle with connection  
  $\theta$. For $P\in I^l_0(G)$,  let $\widetilde{TP(\theta)}$ be the
  mod $\Z$ reduction of the real cochain $TP(\theta)$. Then there exists a
 cochain $U\in C^{2l-1}(B,\R/\Z)$ such that
\begin{equation*}
\widetilde{TP(\theta)}=\pi^*(U)+\text{ coboundary}.
\end{equation*}
\end{thm}

The proof is essentially given in Theorem \ref{thm:CohInf} below.  

\begin{cor}\label{co:CSForm}\cite[Theorem 3.16]{C-S} Assume $P\in I^l_0(G)$
and
  $P(\Omega^l_{E})=0$. Then there exists $CS_P(\theta) \in H^{2l-1}(B,\R/\Z)$
  such that
\begin{equation*}\label{one}
\left[ \widetilde{ TP(\theta)}\right]=\pi^*(CS_P(\theta)).
\end{equation*}
\end{cor}

\begin{proof}
Choose $U\in C^{2l-1}(B,\R/\Z)$ as in Theorem \ref{thm:Coh}. Since
$P(\Omega^l_{E})=0$, Proposition \ref{eq:CS} implies
$\delta\widetilde{TP(\theta)}=\widetilde{ dTP(\theta)}=0$.  By Theorem
\ref{thm:Coh}, $\pi^*U$ and $\widetilde{ TP(\theta)}$ are 
cohomologous. Set $CS_P(\theta)=[U]$.
\end{proof}

Notice that the secondary class or {\it Chern-Simons class} $CS_P(\theta)$, is
defined only when the characteristic form $P(\Omega_E)$
vanishes.  The proof of Theorem \ref{thm:CohInf} shows that $CS_P(\theta)$ is
independent of the choice of $U$.  

The following corollary will be taken as the definition of Chern-Simons
classes for trivial $\G$-bundles (see Definition \ref{de:CS2}).

\begin{cor} \label{trivialbundles} Let $(E,\theta)\to B$ 
be a trivial $G$-bundle with connection, and let $\chi$ be a global section.
For  $P\in I^l_0(G)$,
$$CS_P(\theta) = \chi^*[\widetilde{TP(\theta)}].$$
\end{cor}

\begin{proof}  This follows from Corollary \ref{co:CSForm} and $\pi\chi = {\rm
    Id}.$  \end{proof}

\subsection{Chern-Simons Theory on Loop Spaces}

In \cite{P-R2}, Chern forms are defined on vector bundles with structure group
$\G$ and with $\G$-connections, or equivalently on   principal $\G$-bundles
with connections, where the $\Pdo$s act
on sections of a finite rank bundle $E\to N$ over a closed manifold.
The key technical point is to find suitable analogs for the
polynomials $P\in I^l(G)$.  We single out two 
analogs of the Newton
polynomials $\Tr(A^l)$:
for $A\in \G$, define
\begin{equation}\label{eq:InfTr0} 
P_l^{(0)}(A)=k(l)\int_{S^*N}\Tr
  \left(\sigma_0(A^l)(x,\xi)\right)\ d\xi dx.
\end{equation}
Here $S^*N$ is the unit cosphere bundle of $N$ and $k(l)=(2\pi
i)^{-l}({\rm Vol}\ N)^{-1}.$ 
Note that $d_l = (2\pi i)^{-l}$ is the normalizing constant
such that
$[d_l\tr((\Omega^u)^l)]\in H^{2l}(BU(n),\Z)$
for a connection $\theta^u$ on 
$EU(n)\to BU(n)$.
 In \cite{P-R}, $P_l^{(0)}$ is
called a {\it Leading Order Symbol Trace}.

The second analog is
\begin{equation}\label{eq:InfTr1} 
P_l^{(-1)}(A)=
k(l)
  i^{n}
\int_{S^*N}\Tr \left(\sigma_{-n}(A^l)
(x,\xi)\right)\ d\xi dx. 
\end{equation}
$P_l^{(-1)}(A)$ is a multiple of the {\it Wodzicki residue} of
$A^l$. 
  The factor $i^n$ insures that the Wodzicki residue of a real operator
  is real. 
\medskip
As usual, $P_l^{(i)}, i = 0, -1$,
determine  polynomials by polarization.

For $P_l^{(i)}$, the commutativity and invariance properties
  hold because (\ref{eq:InfTr0}) and (\ref{eq:InfTr1}) are tracial 
  \cite{P-R2} (i.e.~$\Tr[\sigma_i(AB)] = \Tr[\sigma_i(BA)]$ 
  for $A, B \in \Pdo_{\leq 0}$).  In particular, $P_l^{(i)}$ are in
both $I^l(\calG), I^l(\G)$
(although trivially $P_l^{(-1)} = 0$ on $\calG$).

The proof of Proposition
  \ref{eq:CS} to carries over to $\G$-bundles, and so to
  the frame bundle of loop space.

\begin{prop}\label{pr:CSInf}
For a bundle $\G$-bundle with connection
$(\calE,\theta)\overset{\pi}{\to} \calB$,
and for $P\in I^l(\G)$, set
\begin{equation}\label{eq:ChernSimonsInf}
\begin{split}
\phi_t&=t\Omega+\frac{1}{2}(t^2-t)[\theta,\theta],\\
TP(\theta)&=l\int_0^1 P(\theta\wedge \phi^{l-1}_t)dt
\end{split}
\end{equation}
Then $dTP(\theta)=P(\Omega^l)$.
\end{prop}

In the Proposition, we can replace $\G$ by $\calG$.

\begin{rem}  The tracial properties of $P_l^{(i)}$  imply that
  $P_l^{(i)}(\Omega)$ is a closed form with cohomology class
  independent of the connection $\theta.$  
The cohomology classes for
  $P_l^{(0)}, P_l^{(-1)}$ 
are the components of the so-called {\it leading order Chern character} and
the  {\it Wodzicki-Chern character}.  Using Newton's formulas, the Chern
  characters define Chern classes $c_k^{(0)}$, $c_k^{\rm res}$, as usual.
Examples of nontrivial leading
  order Chern classes are given in \cite{P-R2}.  No nonzero examples of
  Wodzicki-Chern classes are known; see \S4.2.
\end{rem}
\medskip

The main goal of this section is to show that 
Theorem \ref{thm:Coh} extends to the
frame bundle $FLM$ for $P = P_l^{(0)}$.  For $P_l^{(-1)}$, we only get an
extension of Corollary \ref{trivialbundles}.

As a first step, we have

\begin{lem}\label{eq:IO}
Let $\G$ be the set of invertible zeroth order $\Pdo$s
  acting on sections of the trivial bundle $\calR^n$. 
Then $P_l^{(0)} \in
I_0(\calG)$.
\end{lem}

\begin{proof}
See \S\ref{se:proof}.
\end{proof}

As in \cite{C-S}, we have

\begin{thm}\label{thm:CohInf} Let $(\calE,\theta)\to \calB$ be a 
$\calG$-bundle with connection
  $\theta$ and assume $P_l(\Omega)=0$. Let
 $\widetilde{ TP(\theta)}$ be the mod $\Z$ reduction of $TP(\theta)$. Then
  there exists a cochain $U\in C^{2l-1}(\calB,\R/\Z)$ such that
\end{thm}
\begin{equation*}
\widetilde{ TP(\theta)}=\pi^*(U)+\text{ coboundary}.
\end{equation*}

\begin{proof}
By \cite[\S 4]{P-R}, 
$E\calG\to B\calG$ has a universal connectionn $\hat\theta$ (with curvature
$\hat\Omega$). Thus there exists a geometric classifying map $\phi: \calB\to
B\calG$: i.e.~$(\calE,\theta)\simeq (\phi^* E\calG,\phi^* \hat\theta)$. By
Lemma \ref{eq:IO}, $P\in I^l_0(\calG)$, so its mod $\Z$ reduction is
zero. From the Bockstein sequence
\begin{equation*}\label{eq:bockstein}
\begin{CD} \cdots\longrightarrow
 H^i(B\calG,\Z)\longrightarrow H^i(B\calG,\R) @>{{\rm mod}\ \Z}>> 
H^i(B\calG,\R/\Z)\longrightarrow H^{i+1}(B\calG,\Z)\longrightarrow\cdots
\end{CD}
\end{equation*}
we deduce that $P(\hat\Omega)$ represents an integral class in $B\calG$. Thus
$\widetilde{P(\hat\Omega)}$ as a cochain vanishes on all cycles in $B\calG$,
and hence is an $\R/\Z$ coboundary, i.e.~there exists
 $\bar u \in C^{2l-1}(B\calG,
\R/\Z)$ such that $\delta \bar u=\widetilde{P(\hat\Omega)}$. We have
\begin{equation*}
\delta \pi^*(\bar u)=\pi^*(\delta u)=\pi^*(\widetilde{P(\hat\Omega)})
= \widetilde{dTP(\hat\theta)}=\delta(\widetilde{TP(\hat\theta)}).
\end{equation*}
The acyclicity of $E \calG$ implies
$\widetilde{TP(\hat\theta)}=\pi^*(\bar u)+\text{coboundary}.$
Now set $U=\phi^*(\bar u).$
\end{proof}

\begin{defn} \label{def:CS1}  
{\it 
Let  $(\calE,\theta)\to \calB$ be a $\calG$-bundle with connection $\theta$
and curvature $\Omega$, and assume
$P_l^{(0)}(\Omega) = 0.$  In the notation of Theorem \ref{thm:CohInf}, define
the }  Chern-Simons
  class $ CS_{2l-1}(\theta)\in H^{2l-1}(B,\R/\Z)$  {\it by  }
$$CS^{(0)}_{2l-1}(\theta) = [U].$$
\end{defn}

If $\G$ acts on $E\to N$, the top order symbol is a homomorphism
$\sigma_0:\G\to\calG$, where $\calG$ acts on $\pi^*E\to S^*N.$  
A $\G$-bundle $\calE$ has an associated $\calG$-bundle $\calE'$
with transition
function $\sigma_0(A)$, if $A$ is a transition function of $\calE$.  A
connection $\theta$ with curvature $\Omega$
on $\calE$ gives rise to a connection $\theta' = \sigma_0(\theta)$
on $\calE'$ with curvature
$\sigma_0(\Omega).$  Since $P_l^{(0)}(\Omega) = P_l^{(0)}(\sigma_0(\Omega))$,
we define $CS_{2l-1}(\theta) = CS_{2l-1}(\theta')$.

This indirect definition is necessary at present, because we 
do not know if $E\G\to B\G$ admits a universal connection.  As a result, we
can only extend the classical definition of Chern-Simons classes to
$P_l^{(-1)}$ for
trivial $\G$-bundles, using the construction of Corollary
\ref{trivialbundles}.  

\begin{defn}\label{de:CS2} 
{\it 
For trivial $\G$-bundles with connection $(\calE,\theta)\to\calB$ and global
section $\chi:\calB\to\calE$,
and assume that $P_l^{(-1)}(\Omega) = 0.$  Then
the}  Chern-Simons
  class $ CS_{2l-1}(\theta,\chi)\in H^{2l-1}(\calB,\R/\Z)$ {\it is defined by }
\begin{equation*}
CS_{2l-1}^{(-1)}(\theta,\chi) = \chi^*\left[ \widetilde{ TP(\theta)}\right].
\end{equation*}
\end{defn}

\begin{rem}
The Chern-Simons class is independent of 
  the section $\chi$ for finite dimensional groups and for $\calG$, since it
  is defined via a universal connection.  
\end{rem}

\section{{\bf Properties of Wodzicki-Chern Classes}}

In this section we give a proof of Lemma \ref{eq:IO}.
We also give a vanishing result for Wodzicki-Chern classes on
mapping spaces of manifolds
generalizing \cite{M}.

\subsection{Integrality of Leading Order Symbol Characteristic Classes}
\label{se:proof}

The goal of this subsection is to show that $W(P_l^{(0)}) \in
H^{2l}(B\calG, \Z)$.

We do not know if the corresponding result $W(P_l^{(-1)})\in H^{2l}(B\G, \Z)$
 is true.
Fortunately, for 
$FLM$, we know that $P_l^{(-1)}(\Omega^s) = 0$  (Lemma
\ref{quickref}).  
  This suffices to define the Chern-Simons class for
the Levi-Civita connections on $FLM$, if $FLM$ is trivial, although the class
depends on the choice of
global section.


By \cite{A-B1},
$B\calG= C^\infty_{(0)}(S^1,BSO(n))=\{f:S^1\to BSO(n) | f^*ESO(n) \simeq
\pi^*\calR^n\}$. 
As a more general setup, consider a closed manifold $N$ and a
finite rank real bundle $E\to M$.
Let
$\ev:C^\infty(N,M) \times N\to M$ be the evaluation map 
$\ev(f,n) = f(n).$  The bundle 
$\ev^*E$ determines an 
infinite rank bundle $\pi_*\ev^*E \to C^\infty(N,M)$, where $\pi_*\ev^*E|_f =
\Gamma(f^*E\to N),$ with $\Gamma$ denoting some Sobolev space of sections.  
(Here $\pi:C^\infty(N,M) \times N\to 
C^\infty(N,M)$ is the projection.)
For $n\in N$, define $\ev_{n}: C^\infty(N,M)\to M$ by $\ev_n(f) = f(n).$ 

It is well known that connections push down under $\pi_*$.  For the gauge
group case, this gives the following:

\begin{lem} \label{lem:ConnEG}
The universal bundle $E\calG\to B\calG$ is isomorphic to $\pi_*\ev^*ESO(n)$.
$E\calG$ has a universal connection
$\theta^{E\calG}$ defined on $s\in\Gamma(E\calG)$ by
\begin{equation*}
(\theta^{E\calG}_Z s)(\gamma)
(\alpha)=\left( (\ev^*\theta^u)_{(Z,0)} u_s\right)(\gamma,\alpha).
\end{equation*}
Here $\theta^u$ is the universal connection on $ESO(n)\to BSO(n)$, and
$u_s: C^\infty(N,M)\times N \to \ev^*ESO(n)$ is defined by
$u_s(f,n)=s(f)(n)$.
\end{lem}

\begin{proof}
See \cite{P-R}.
\end{proof}

\begin{cor}
The curvature $\Omega^{E\calG}$ of $\theta^{E\calG}$ satisfies
\begin{equation*}\label{eq:pullback}
\Omega^{E\calG}(Z,W)s(f)(n)=\ev^*\Omega^u  ((Z,0),(W,0)) u_s(f,n).
\end{equation*}
\end{cor}

\begin{proof}
This follows from
\begin{equation*}
\Omega^{E\calG}(Z,W)s(f)(n)= [\nabla^{E\calG}_Z \nabla^{E\calG}_W
-\nabla^{E\calG}_W \nabla^{E\calG}_Z -\nabla^{E\calG}_{[Z,W]}]
s(f)(n)
\end{equation*}
and the previous lemma.
\end{proof}

We now prove that $P_l^{(0)}(\Omega^{E\calG})\in H^{2l}(B\calG,\Z)$.

Since
${\rm ev}_{0} = \ev_{\alpha_0}$ is homotopy equivalent for every $\alpha_0\in
S^1$, the cohomology class
\begin{equation*}\label{eq:Inde}
\left[ P_l^0(\ev_0^*\Omega^u)\right]\in H^{2k}(B\calG\times\{n_0\},\R)
 \cong H^{2k}(B\calG,\R)
\end{equation*}
is
independent of $\alpha_0$.  Thus
\begin{eqnarray}\label{eq:IntPO}
\left[\frac{d_l}{4\pi}\int_{S^*S^1} \tr\ \sigma_0((\Omega^{E\calG})^l) d\xi
d\alpha \right] &=& \frac{d_l}{4\pi}\int_{S^*S^1} \left[\tr
\ \sigma_0 \left((\ev_{\alpha}^*\Omega^u)^l\right)
\right] d\xi d\alpha, \nonumber \\
&=& \left[d_l \ev^*_{0}\tr\  \sigma_0
\left((\Omega^u)^l\right) \right], \\ 
&=& \ev^*_{0}\left[d_l \tr
\left((\Omega^u)^l\right)\right], \nonumber   
\end{eqnarray}
since $\Omega^u$ is a multiplication operator.  By the choice of $d_l$, 
the last term in (\ref{eq:IntPO}) lies in
$\ev^*_{0} H^{2l}(BSO(n),\Z)
\subset H^{2l}(B\calG),\Z)$.
Thus
\begin{equation*}\label{eq:WP0}
 W(P_l^{(0)}) = [P_l^{(0)}(\Omega^{E\calG})]\in H^{2l}(B\calG,\Z),
\end{equation*}
which completes the proof of the Lemma \ref{eq:IO}.\\

\medskip
The following table summarizes the results.

\bigskip

\centerline{
\begin{tabular}{c|c @{}*{3}{c}}
\multicolumn{3}{c}{{\sc Table} 1. Is\ $P_l^{(i)}\in I_0(G)$?} \\
\multicolumn{5}{c}{ } \\
$\downarrow P_l^{(i)}\ \ \vline \ \ G\rightarrow$  & $\calG$ &&& $\G$ \\
\hline
&&&&\\
$P_l^{(0)}$  & yes&& & ?\\
&&&&\\
$P_l^{(-1)}$ &yes, trivially &&& ? \\
&&&& (but see Lemma 4.5) \\
\end{tabular}
}

\bigskip

\begin{rem}  Let $(\calE,\theta)\to\calB$ be a $\calG$-bundle with connection,
  where $\calG$ is the gauge group of the rank $n$
hermitian bundle $E\to N$, 
  and let $f:\calB\to B\calG$ be a geometric classifying map.  The argument
  above easily extends to show that the $l^{\rm th}$ leading order Chern class
  equals $f^*\ev_0^*c_l(EU(n)).$  Thus all leading order Chern classes are
  pullbacks of finite dimensional Chern classes.  (This argument was developed
  with S. Paycha.)
\end{rem}

\subsection{A Vanishing Theorem for Wodzicki-Chern Classes}

\begin{thm}  
If  $\calE\to C^\infty(N,M)$ satisfies ${\calE} =
\pi_* ev^*E$ as above, then the Wodzicki-Chern classes $c_k^{res}({\calE})$
vanish for all $k$.  
\end{thm}

\begin{proof}
As in the previous subsection, $\calE$ admits a connection whose curvature
$\Omega$ is a multiplication operator.  $\Omega^l$ is also a multiplication
operator, and hence $P_l^{(-1)}(\Omega) = 0.$
\end{proof}

For a real infinite rank bundle, Wodzicki-Pontrjagin classes are defined as in
finite dimensions:  $p_k^{\rm res}(\calE) = (-1)^k c_{2k}^{\rm
  res}(\calE\otimes \C).$  

\medskip
\begin{cor}  The Wodzicki-Pontrjagin classes of $TC^\infty(N,M)$ and of all
  naturally associated bundles vanish.
\end{cor}

\begin{proof}  Pick an element $f_0$ in a fixed path component $A_0$ of
  $C^\infty(N,M).$  For $f\in A_0$, $T_fC^\infty(N,M) \simeq \Gamma(f^*TM\to
  N)\simeq \Gamma(f_0^*TM\to   N)$, where the second isomorphism is
  noncanonical.  Thus over each component, $TC^\infty(N,M)$ is of the form
  $\pi_*\ev^*TM.$  The vanishing of the Wodzicki-Pontrjagin classes of
  associated bundles (such as exterior powers of the tangent bundle) follows
  as infinite dimensions, since there is a universal geometric bundle.
\end{proof}
\medskip

We also have a trivial vanishing result for Wodzicki-Pontrjagin forms for
$FLM.$  
\medskip

\begin{lem} \label{quickref}
 The forms $P_l^{(-1)}(\Omega^s)$, $l>1$, vanish on $FLM.$ 
\end{lem}

\begin{proof}  This follows from the fact that $(\Omega^s)^l$  is of order
  $-l.$  
\end{proof}

\begin{rem} Similarly, if $\calE\to\calB$ is an infinite rank $\G$-bundle, for
  $\G$ acting on
$E\to N^n$, and if $\calE$ admits a 
   $\G$-connection whose
  curvature has order $-k$, then $0 =
  c_{[n/k]}(\calE) =   c_{[n/k]+1}(\calE) = \ldots$
  Thus the Wodzicki-Chern classes are obstructions to
  the negativity of the order of the curvature.
\end{rem}

\section{{\bf The Chern-Simons Class and Parallelizable Manifolds}}

Thanks to Chern-Simons formalism, 
the vanishing of the curvature expressions $P_l^{(-1)}(\Omega^s)$, $l>1$,
in the last lemma is actually an advantange:
as in \S3, we can define 
Chern-Simons $CS_{2l-1}(\theta^s,\chi)\in H^{2l-1}(LM,\R/\Z)$ provided $LM$ is
parallelizable.

In this section, we
describe a strategy to detect non-trivial Chern-Simons classes
on parallelizable loop spaces. 

\begin{lem}\label{le:41}
If $M$ is parallelizable, then $LM$ is  parallelizable.
\end{lem}

\begin{proof}
Let $\phi: TM\to M\times \R^n$ be a trivialization of
$TM$.  For
$X_\gamma\in T_\gamma LM=\Gamma(\gamma^*TM)$, define
\begin{equation*}
\begin{split}
\Psi :TLM&\to LM\times \Gamma(S^1\times\R^n\to S^1)\\
X_\gamma &\longmapsto (\gamma, \alpha\mapsto \pi_2(\phi(X_\gamma(\alpha)))),
\end{split}
\end{equation*}
where $\pi_2:M\times \R^n\to \R^n$ is the projection. 
It is easy to check that $\alpha$ is a smooth  trivialization of $TLM$ in the
$H^s$ norm.
\end{proof}

Therefore, for parallelizable $M$ there exists a global section $\chi:LM\to
FLM$. For $P$ equal either $P_l^{(0)}$ or $P_l^{(-1)}$ and $l>1$,
$\chi^*TP(\theta^s)\in H^{2l-1}(LM,\R)$ and $[\widetilde{ \chi^*
TP(\theta^s)}] = CS^{(i)}_{2l-1}(\theta^s,\chi)$ for $i = 0, 1.$
Thus $CS^{(i)}_{2l-1}(\theta^s) = 
CS^{(i)}_{2l-1}(\theta^s,\chi)$  
is nontrivial if
there exists\\
 $[z]\in H_{2l-1}(LM;\Z)$ with 
\begin{equation*}\langle
\chi^*TP(\theta^s),[z] \rangle\not\in \Z.\end{equation*} 
To find an
appropriate $[z]$, consider the map
\begin{equation*}
\beta: N\to L(N\times S^1),\ \ 
x\mapsto \left(\beta (x)(\alpha)=(x,\alpha)\right).
\end{equation*}

\begin{lem}\label{eq:inj}
$\beta_*:H_i(N,\Z)\to H_i(L(N\times S^1),\Z)$ is injective.
\end{lem}

\begin{proof}
Fix $\alpha_0\in S^1$ with its associated evaluation map $ \ev_0 =
\operatorname{ev}_{\alpha_0}: L(N\times S^1)\to N\times S^1$.
Let $\pi_1:N\times S^1\to N$ be the projection. From 
\begin{equation*}
\begin{CD}
N@>\beta>>L(N\times S^1)@>{\rm ev}_{0}>> N\times S^1@>\pi_1>> N,
\end{CD}
\end{equation*}
we obtain $
\pi_1 \circ {\rm ev}_{0}\circ \beta={\rm Id}_N$, which
implies that $\beta_*$ is injective.
\end{proof}

Set
$N=S^3$ and $M = S^3\times S^1.$
By the lemma,  $\beta_*[S^3]\in H_3(L(S^3\times S^1),\Z)$ is
nontrivial. This class works well for $l=2$, since
$CS^{(i)}_3(\theta^s)\in H^3(LM,\R/\Z)$.

\begin{cor}\label{bigcor} $CS^{(i)}_3(\theta^s)$ is
nontrivial in $H^3(L(S^3\times S^1),\R/\Z)$ if
\begin{equation}\label{eq:IntS3}
\begin{split}
\langle \chi^*TP(\theta),\beta_*[S^3] \rangle= \int_{S^3} \beta^*\chi^*
TP(\theta)\not \in \Z.
\end{split}
\end{equation}
\end{cor}

\begin{rem}
To compute the
integrand in (\ref{eq:IntS3}), it is useful to pick a global frame
$\{E_1,E_2,E_3\}$ of $S^3$ and $E_4 = \partial_\alpha$ for $S^1.$ Then of
course 
\begin{equation*}\label{eq:pull}
\beta^*\chi^*TP(\theta)(E_1,E_2,E_3)=
\chi^*TP(\theta)(\beta_*E_1,\beta_*E_2,\beta_*E_3).
\end{equation*}
It is easy to check that
\begin{equation}\label{eq:frame}
\beta_*(E_1)
= ( E_1,0,0,0), \beta_*(E_2)=(0, E_2,0,0,0),\ \beta_*(E_3)=(0,0,E_3,0)
\end{equation}
as constant sections of the trivial bundle $T_{\beta(m)}L(S^3\times S^1)$.
\end{rem}
\medskip
${}$

\section{{\bf Calculations on $S^3\times S^1$.}}

In this section we explicitly compute a three dimensional Chern-Simons class.
In \S\ref{sec1}, we begin the computations, and
show the vanishing of
 the Chern-Simons class $CS_3^{(0)}(\theta^s)\in
H^3(LM,\R/\Z)$ associated to the Levi-Civita connection $\theta^s$ for a class
of metrics
on $M = S^3\times S^1.$  
In \S\ref{sec2}, we find a metric on $M$ and a global frame $\chi$ of $FLM$
such that $CS_3^{(-1)}(\theta^s,\chi) \neq 0$.

\subsection{Computations on $S^3\times S^1$ for  $P_2^{(0)}$}\label{sec1} 

Setting $l=2$ and combining (\ref{eq:InfTr0}) with (\ref{eq:ChernSimonsInf})
  for $\theta=\theta^s$ on the frame bundle $FLM$, we obtain
\begin{equation*}\label{eq:lval}
\begin{split}
\chi^*TP(\theta)&=-\frac{d_2}{6\cdot 4\pi} P^{(0)}(\chi^*\theta
\wedge \chi^*\theta \wedge \chi^*\theta)\\ 
&= -\frac{d_2}{6\cdot 4\pi}
\int_{S^*S^1} \tr \ \left[ \sigma_0(\chi^*\theta) \wedge
\sigma_0(\chi^*\theta) \wedge \sigma_0(\chi^*\theta)\right] d\xi d\alpha.
\end{split}
\end{equation*}
To simplify notation, set $\omega^s=\chi^*\theta$. Once we choose a
metric for $M = S^3\times S^1$ and $\chi$, we can use (\ref{localexp}) and
(\ref{eq:zero}) to explicitly compute
$\chi^*TP(\theta)$. 

Let
$ E_1,E_2, E_3$ be a frame of orthonormal left invariant vector fields 
for $S^3$ with the standard metric. 
Let $E_4=\partial_\rho$, where $\rho$ is the coordinate on $S^1$ in $S^3\times
S^1$,  
and impose the usual Lie relations
\begin{equation*}
[E_1,E_2]=2E_3,\ [E_2,E_3]=2E_1,\ [E_1,E_3]=-2E_2, [E_i,E_4]=0.
\end{equation*}
$\chi:LM\to FLM$ will be the ``loopification'' of the global frame
$(E_i)$ of $S^3\times S^1$: $\chi(\gamma)(\alpha) =
(E_i(\gamma(\alpha)))$, which is identified with a gauge transformation of
$T_\gamma LM$ under the isomorphism $T_\gamma LM \simeq \calR^4.$  

Fix functions $\lambda=\lambda(\alpha)$, $\mu=\mu(\alpha)$,
$\nu=\nu(\alpha)$. Take the metric on $M$
for which $\lambda E_1, \mu E_2, \nu E_3,
E_4$ are orthonormal. The non-zero Christoffel coefficients are
\begin{equation*}
\begin{split}
\Gamma^3_{12}&=\left(\frac{\mu^2\lambda^2-\mu^2\nu^2+\nu^2\lambda^2}
{\lambda\mu\nu}\right)=-\Gamma^2_{13},\
\Gamma^3_{21}=\left(\frac{-\mu^2\lambda^2-\mu^2\nu^2+\nu^2\lambda^2}
{\lambda\mu\nu}\right)=-\Gamma^1_{23},\\
\Gamma^2_{31}&=\left(\frac{\nu^2\lambda^2-\lambda^2\mu^2+\mu^2\nu^2}
{\lambda\mu\nu}\right)=-\Gamma^1_{32}, \
\Gamma^1_{14}= - \Gamma^4_{11}=-\frac{\dot \lambda}{\lambda},\ 
\Gamma^2_{24}= - \Gamma^4_{22}=-\frac{\dot \mu}{\mu},\\
\Gamma^2_{24}&= - \Gamma^4_{22}=-\frac{\dot \nu}{\nu},\ 
\Gamma^i_{4j}=0=\Gamma^i_{44}=\Gamma^4_{4j}=0.
\end{split}
\end{equation*}
 Set
\begin{equation}\label{eq:Coefficients}
\begin{split}
U&=\nu^2\frac{\mu^2-\lambda^2}{\lambda\mu\nu}, \
V=\mu^2\frac{\mu^2-\lambda^2}{\lambda\mu\nu}, \
W=\lambda^2\frac{\nu^2-\mu^2}{\lambda\mu\nu},\\
A&=\frac{\dot \lambda}{\lambda}, 
\ B=\frac{\dot \mu}{\mu}, \ C=\frac{\dot \nu}{\nu}.
\end{split}
\end{equation}
A direct calculation gives 
\begin{equation*}\label{eq:MaTr0}
\sigma_0 (\omega^s) = 
\left(\begin{array}{rrrr}-A\psi^4& U\psi^3 &-V \psi^2& \frac{1}{2}A\psi^1\\ 
U\psi^3&-B\psi^4 &W \psi^1& \frac{1}{2}B\psi^2\\ 
-V\psi^2&W\psi^1&-C\psi^4&\frac{1}{2}C\psi^3 \\
\frac{1}{2}A\psi^1&\frac{1}{2}B\psi^2&\frac{1}{2}C\psi^3 &
0\ \ \end{array}\right)
\end{equation*}
Here $\{\psi^i\}$ is the frame dual  to $\{E_i\}$ 


A straightforward computation using (\ref{eq:zero1}) gives
$\Tr\left[\sigma_0(\omega^s)\wedge \sigma_0(\omega^s)
\wedge\sigma_0(\omega^s)\right] = 0.$
Thus $CS_3^{(0)}(\theta^s) = 0$ for this class of metrics, so we turn our
attention to $CS_3^{(-1)}(\theta^s).$  

\subsection{Computations on $S^3\times S^1$ for $P^{(-1)}$}\label{sec2}

For the case $l=2$, (\ref{eq:ChernSimonsInf}) gives
\begin{eqnarray}\label{eq:CSP1}
\chi^*TP(\theta)&=&
2 \int_0^1 P_2^{(-1)}\left(\chi^*\theta\wedge t\chi^*\Omega 
+\frac{1}{2}(t^2-t)[\chi^*\theta,\chi^*\theta]\right) dt\nonumber\\
&=& P_2^{(-1)}(\chi^*\theta\wedge \chi^*\Omega)
-\frac{1}{6} P^{(-1)}(\chi^*\theta\wedge\chi^*\theta\wedge\chi^*\theta)\\
& =& -\frac{i}{8\pi^3}
\int_{S^*S^1}\Tr \left[\sigma_{-1}
(\chi^*\theta\wedge \chi^*\Omega)\right] d\xi d\alpha\nonumber\\ 
&&\ \ +\frac{i}{48\pi^3 }
\int_{S^*S^1}\Tr \left[\sigma_{-1}
(\chi^*\theta\wedge \chi^*\theta\wedge\chi^*\theta)\right]  d\xi d\alpha,
\nonumber
\end{eqnarray}
where $\theta=\theta^s, \Omega=\Omega^s$.
By the symbol calculus for $\Pdo$s and (\ref{eq:zero1}), we have 
\begin{equation}\label{eq:Conn}
\Tr \left[\sigma_{-1}(\chi^*\theta\wedge\chi^*\theta\wedge\chi^*\theta)\right]
= 3\Tr \left[\sigma_{-1}(\chi^*\theta)\wedge\sigma_0(\chi^*\theta)\wedge 
\sigma_0(\chi^*\theta)\right]
\end{equation} 
\begin{equation}\label{eq:Curv}\sigma_{-1}(\chi^*\theta\wedge\chi^*\Omega)
=\sigma_0(\chi^*\theta)\wedge\sigma_{-1}(\chi^*\Omega)
\end{equation}
As in \S\ref{sec1}, we may replace $\chi^*\theta$ by $\omega^s$ and
$\chi^*\Omega$ by $\Omega^s.$  

First, we compute the contribution from (\ref{eq:Conn}) to (\ref{eq:CSP1}).
Recall that we need to compute the terms on  the right hand side of
 (\ref{eq:Conn})  on $\beta_*(TS^3)$. 
On the loop $\beta(m)(\alpha)=(m,\alpha)$, we have
$\partial_\alpha \gamma=\dot \gamma =(0,0,0,1)$. Thus (\ref{eq:Minus}) becomes
%
\begin{eqnarray}\label{eq:Minus1}
\sigma_{-1}(\omega^s(X))^a_b &=& 
2is \xi^{-1}\left[ \partial_l
\Gamma^a_{b4}+\left(\Gamma^a_{lk}\Gamma^k_{b4}
-\Gamma^a_{k4}\Gamma^k_{lb}
-\delta^{ap}\delta_{br}\Gamma^r_{q4}\Gamma^q_{pl}
-\delta^{pq}\delta_{rb}\Gamma^a_{p4}\Gamma^r_{ql}\right)\right]X^l\nonumber\\
&&\ +2is\xi^{-1}\left[(\Gamma^a_{bl}+\delta^{ap}\delta_{qb}
\Gamma^q_{pl})\partial_\alpha X^l
  +\partial_\alpha\left(\Gamma^a_{lb}+\delta^{ap}
\delta_{qb}\Gamma^q_{pl}\right)X^l\right].
\end{eqnarray}
Note that $\partial_\alpha X^l=0$ for $X\in  \beta_*(TS^3)$, 
as $\beta_*(E_i)$ does not
depend on $\alpha$
by (\ref{eq:frame}).
Thus on $\beta_*(TS^3)$, (\ref{eq:Minus1}) reduces to
\begin{eqnarray}
\label{eq:AB}
\sigma_{-1}(\omega^s(X))^a_b &=& 
2is \xi^{-1}\left[ \partial_l
\Gamma^a_{b4}+\left(\Gamma^a_{lk}\Gamma^k_{b4}-
\Gamma^a_{k4}\Gamma^k_{lb}
-\delta^{ap}\delta_{br}\Gamma^r_{q4}\Gamma^q_{pl}
-\delta^{pq}\delta_{rb}\Gamma^a_{p4}\Gamma^r_{ql}\right)\right]X^l\nonumber\\
&&\ +2is\xi^{-1}\left[
  \partial_\alpha\left(\Gamma^a_{lb}+\delta^{ap}
\delta_{qb}\Gamma^q_{pl}\right)X^l\right]
\end{eqnarray}
Combining (\ref{eq:AB}) 
with the values of the
Christoffel symbols from \S6.1 gives a messy but explicit expression for the
contribution of (\ref{eq:Conn}) to (\ref{eq:CSP1}).

Second, we compute the contribution of (\ref{eq:Curv}) to
(\ref{eq:CSP1}). The $-1$ order symbol of the
curvature in our orthonormal frame is 
\begin{eqnarray}\label{eq:Curv1}
\frac{1}{2is \xi^{-1}}\sigma_{-1}(\Omega^s(X,Y))^k_l &=& 
\dot X^p Y^r\left[
  \partial_p\Gamma^k_{r l} -\partial_r\Gamma^k_{p
  l}-\delta^{kb}\delta_{ml}\partial_r\Gamma^m_{b p}\right]\nonumber\\
&& \ \ + X^pY^r
  \left[\partial_{p4}\Gamma^l_{r m}+\delta^{kb}\delta_{ml}
\partial_{p4}\Gamma^m_{b r}-
  \partial_{r 4}\Gamma^k_{pl}-
\delta^{kb}\delta_{ml} \partial_{r 4}\Gamma^m_{bp}\right]\nonumber\\ 
&& \ \  +X^p \dot Y^r\left[ \partial_p\Gamma^k_{l r} 
+ \delta^{kb}\delta_{ml} \partial_p\Gamma^m_{b  r}
-\partial_r\Gamma^k_{pl}\right]
\end{eqnarray}
\cite{M-R}. Arguing as above (\ref{eq:AB}), we see
that the first and third terms on the right hand side of 
(\ref{eq:Curv1}) do not contribute to the Chern-Simons
class. In the second term, the only possible
contributions come from $X^k = X^4$ or $Y^r = Y^4$, but again by
(\ref{eq:frame}), 
these  components vanish on $ \beta_*(TS^3)$. 
Thus  (\ref{eq:Curv}) does not contribute to
(\ref{eq:CSP1}).

In
summary, on this image
\begin{equation*}\label{eq:FinInt}
\chi^*TP(\theta)=
\frac{s}{2\pi^2}\psi_1\wedge\psi_2\wedge\psi_3 \int_{S^1}
f\left(\lambda(\alpha),\mu(\alpha), \nu(\alpha)\right),
\end{equation*}
where the complicated function $f$ is determined explicitly by 
(\ref{eq:CSP1}), (\ref{eq:Conn}), (\ref{eq:AB}).
For
\begin{equation*}\label{eq:fam}
\lambda(\alpha)=1, \ \mu(\alpha)=2+ \frac{1}{a}\cos(a\alpha)\sin(a\alpha), 
\  \nu(\alpha)=2-\cos(a\alpha),\  a\in\Z,
\end{equation*}
we compute via Mathematica \cite{Bo} that 
\begin{equation*}\label{eq:FinCom}
\int_{S^3} \beta^*\chi^*
TP(\theta)= \int_{S^3} \frac{s}{2\pi^2}
\psi_1\wedge\psi_2\wedge\psi_3 \int_{S^1}
f(\lambda, \mu,\nu)=\frac{s}{4}\int_{S^1}
f(\lambda, \mu, \nu)\not\in \Z
\end{equation*}
for various choices of $a$. (See Figures 1, 2.)
By Corollary \ref{bigcor},
for these (and many other) choices of $\lambda, \mu, \nu$,
\begin{equation*}
CS_{P^{(-1)}}(\theta)\in H^3(L(S^3\times S^1),\R/\Z)
\end{equation*}
is nontrivial.\newline

\begin{figure}
\epsfig{file=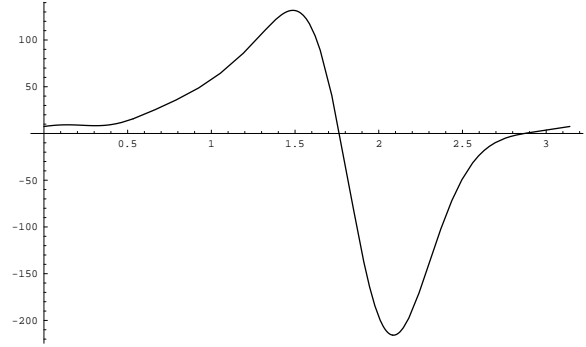,width=3in}
\caption{$a=2$, $\int_{S^1} f(\alpha)=-26.0687$}
\end{figure}

\begin{figure}
\epsfig{file = 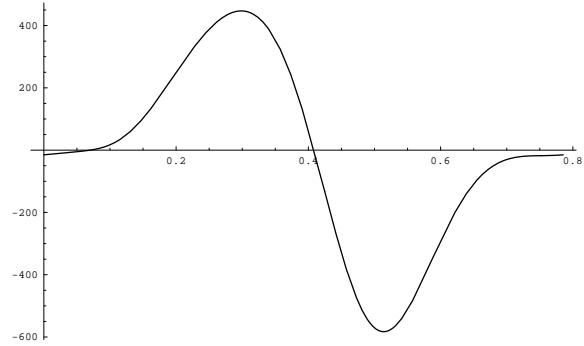,width=3in}
\caption{$a=8$, $\int_{S^1} f(\alpha)=-100.992$}
\end{figure}

\begin{rem}  The Chern-Simons class $CS_3^{(-1)}(\theta,\chi)$
has a linear dependence on the Sobolev
  parameter $s$, and so keeps track of
the $s$-dependence of the
  topology of the frame bundle.  Alternatively, one can
  define 
  a regularization/parameter independent
 Chern-Simons form as $\frac{1}{s}CS_3^{(-1)}(\theta,\chi)$ and note
  that this invariant is non-zero in our example.
\end{rem}

\bibliographystyle{amsplain}
\bibliography{Paper}

\bigskip
\hfill \today \\
\end{document}